\documentclass[11pt]{article}
\usepackage{amsmath,amsthm}

\newtheorem{theorem}{Theorem}
\newtheorem*{lemma}{Lemma}

\title{Monochromatic Progressions in Random Colorings}

\author{Sujith Vijay\\
\small School of Mathematics\\[-0.4ex]
\small Indian Institute of Science Education and Research \\[-0.4ex]
\small Thiruvananthapuram-695016, Kerala, India.\\[-0.4ex]
\small \texttt{sujith@iisertvm.ac.in}}
\date{}

\begin{document}
\maketitle

\begin{center}
{\bf{Abstract}}
\vskip 10pt
\end{center}

Let $N^{+}(k)= 2^{k/2} k^{3/2} f(k)$ and $N^{-}(k)= 2^{k/2} k^{1/2} \,
g(k)$ where $f(k) \rightarrow \infty$ and $g(k) \rightarrow 0$
arbitrarily slowly as $k \rightarrow \infty$. We show that the probability
of a random $2$-coloring of $\{1,2,\ldots,N^{+}(k)\}$ containing a
monochromatic $k$-term arithmetic progression approaches $1$, and the
probability of a random $2$-coloring of $\{1,2,\ldots,N^{-}(k)\}$
containing a monochromatic $k$-term arithmetic progression approaches $0$,
as $k \rightarrow \infty$. This improves an upper bound due to Brown
\cite{Bro99}, who had established an analogous result for $N^{+}(k)= 2^k
\log k f(k)$.

\section*{Introduction}

One of the earliest results in Ramsey theory is the theorem of van der
Waerden \cite{Wae27}, stating that for any positive integer $k$, there
exists an integer $W(k)$ such that any $2$-coloring of $\{1,2,\ldots,
W(k)\}$ yields a monochromatic $k$-term arithemtic progression. The exact
values of $W(k)$ are known only for $k \le 6$. Berlekamp \cite{Ber68}
showed that $W(p+1) \ge p 2^p$ whenever $p$ is prime, and Gowers
\cite{Gow01} showed that $W(k)$ is bounded above by a tower of finite
height, i.e., $$W(k) \le 2^{2^{2^{2^{2^{k+9}}}}}$$

Since the best known upper and lower bounds on $W(k)$ are far apart, a lot
of work has been done on variants of the original problem. A natural
question from a probabilistic perspective is to obtain upper bounds on the
slowest growing function $N^{+}(k)$ such that the probability of a
$2$-coloring of $\{1,2,\ldots, N^{+}(k)\}$ containing a monochromatic
$k$-term arithmetic progression (hereafter abbreviated as $k$-AP)
approaches $1$ as $k \rightarrow \infty$. Similarly, one could seek lower
bounds on the fastest growing function $N^{-}(k)$ such that the
probability of a $2$-coloring of $\{1,2,\ldots, N^{-}(k)\}$ containing a
monochromatic $k$-AP approaches $0$ as $k \rightarrow \infty$. An upper
bound for $N^{+}(k)$ was established by Brown \cite{Bro99} who showed that
$N^{+}(k) = O(2^k \log k f(k))$ where $f(k) \rightarrow \infty$ as $k
\rightarrow \infty$. We improve this bound to $N^{+}(k) = O(2^{k/2}
k^{3/2} f(k))$, and also show that $N^{-}(k) = \Omega(2^{k/2} k^{1/2} \,
g(k))$ where $g(k) \rightarrow 0$ as $k \rightarrow \infty$.

\section*{Almost Disjoint Progressions}

A family of sets ${\cal{F}}=\{S_1,S_2,\ldots,S_m\}$ is said to be {\em
{almost disjoint}} if any two distinct elements of ${\cal{F}}$ have at
most one element in common, i.e., if $|S_i \cap S_j| \le 1$ whenever $i
\neq j$. \\

\begin{lemma}
\label{lem1} 
Let ${\cal{F}}_{k,n}$ be the collection of
$k$-APs contained in $\{1,2,\ldots,n\}$ with common difference $d$
satisfying $n/k \le d < n/(k-1)$. Then ${\cal{F}}_{k,n}$ is an almost
disjoint family. Moreover, $|{\cal{F}}_{k,n}| = n^2(1+o(1))/2k^3$.
\end{lemma}

\begin{proof}
Let $n/k \le d < n/(k-1)$, and let
$A_1=\{a,a+d,\ldots,a+(k-1)d\}$ be a $k$-term arithmetic progression in
${\cal{F}}_{k,n}$. For $0 \le \ell \le k-1$, consider the pairwise
disjoint half-open intervals $I_{\ell}=(\ell n/k, (\ell+1)n/k]$. We claim 
that $a+ \ell d \in I_{\ell}$ for $0 \le \ell \le k-1$. Clearly, $a + \ell 
d > \ell d \ge \ell n/k$. Moreover, $a + \ell d = a+(k-1)d - (k- \ell -1)d 
\le n - (k- \ell -1)n/k = (\ell + 1)n/k$. In particular, $a \le n/k \le 
d$. \\

Now suppose that $A_2=\{a',a'+d',\ldots,a'+(k-1)d'\} \in {\cal{F}}_{k,n}$
with $|A_1 \cap A_2| \ge 2$. Let $a + \ell_1 d = a' + \ell'_1 d'$ and $a +
\ell_2 d = a' + \ell'_2 d'$. Since the intervals $I_{\ell}$ are pairwise
disjoint, it follows that $\ell_1 = \ell'_1$ and $\ell_2 = \ell'_2$. But
then we have $a=a'$ and $d=d'$, so that $A_1 = A_2$. Thus
${\cal{F}}_{k,n}$ is an almost disjoint family. \\

Finally, 
$$|{\cal{F}}_{k,n}| = \sum_{\frac{n}{k} \le d < \frac{n}{k-1}} (n-d(k-1))  
= \frac{n^2(1+o(1))}{2k^3}, $$
since there are $n-d(k-1)$ $k$-term arithmetic progressions of
common difference $d$ completely contained in $\{1,2,\ldots,n\}$.
\end{proof}

For each integer $k \ge 3$, let $c_k$ denote the asymptotic constant such
that the size of the largest family of almost disjoint $k$-term arithmetic
progressions contained in $[1,n]$ is $c_k n^2/(2k-2)$. It follows from the
above lemma that $c_k \ge 1/k^2$. Perhaps there is an absolute constant
$\lambda$ such that $c_k \le \lambda/k^2$. Ardal, Brown and Pleasants
\cite{Ard05} have shown that $0.476 \le c_3 \le 0.485$.

\section*{Monochromaticity: Almost Surely and Almost Never}

\begin{theorem}
Let $N^{+}(k)= 2^{k/2} k^{3/2} f(k)$ where 
$f(k) \rightarrow \infty$ arbitrarily slowly as $k \rightarrow
\infty$. Then the probability that a $2$-coloring of
$\{1,2,\ldots,N^{+}(k)\}$ chosen randomly and uniformly contains a
monochromatic $k$-term arithmetic progression approaches $1$ as $k
\rightarrow \infty$.
\end{theorem}

\begin{proof}
Our approach will be similar to that of
Brown\cite{Bro99}, but rather than work with a family of combinatorial
lines in a suitably chosen hypercube, which is an almost disjoint family
of size $O(n)$, we work with $k$-APs of large common difference, which is
an almost disjoint family of size $\Omega(n^2/k^3)$, as shown in the
previous section. \\

Let $n= N^{+}(k)=2^{k/2} k^{3/2} f(k)$ and $q=\lfloor (f(k))^{4/3}
\rfloor$. Let $s=s(k)$ satisfy $n = qs + r, \, 0 \le r < s$. We divide the
interval $[1,n]$ into $q$ blocks $B_1,B_2,\ldots,B_{q}$ of length $s$,
and possibly one residual block $B_{q+1}$ of length $r$.  Let ${\cal{F}}_1
= {\cal{F}}_{k,s}$ consist of all $k$-APs in $B_1=[1,s]$ with common
difference $d$ satisfying $s/k \le d < s/(k-1)$. By Lemma~\ref{lem1}, the
elements of ${\cal{F}}_1$ are almost disjoint, and $s^2/4k^3 \le
|{\cal{F}}_1| \le s^2/k^3$ for large $k$. \\

For each arithmetic progression $P \in {\cal{F}}_1$, let $C_P$ denote the
set of $2$-colorings of $B_1$ in which $P$ is monochromatic. Then
$|C_P|=2^{s-k+1}$. Also, $|C_P \cap C_Q| = 2^{s-2k+2}$, since $|P \cap Q|
\in \{0,1\}$. By Bonferroni's inequality, $$\bigg| \bigcup_{P \in
{\cal{F}}_1} C_P \bigg| \ge \sum_{P \in {\cal{F}}_1} |C_P| -
\sum_{\stackrel{P, Q \in {\cal{F}}_1}{P \neq Q}} |C_P \cap C_Q| =
|{\cal{F}}_1| 2^{s-k+1} - {|{\cal{F}}_1| \choose 2} 2^{s-2k+2}$$

Since $$\frac{|{\cal{F}}_1|}{2^k} \le \frac{n^2}{q^2 2^k k^3} 
\rightarrow 0 \mbox{ as } k \rightarrow \infty , $$ it follows that 
$$\bigg|\bigcup_{P \in
{\cal{F}}_1} C_P \bigg| > |{\cal{F}}_1| 2^{s-k} \ge \frac{2^s
s^2}{2^{k+2} k^3}.$$ Similarly, we can consider the blocks $B_2, B_3,
\ldots, B_q$ and the corresponding families ${\cal{F}}_2, {\cal{F}}_3,
\ldots, {\cal{F}}_q$. Let $p_0$ be the probability that no arithmetic
progression from any of the ${\cal{F}}_i$ is monochromatic under a
$2$-coloring chosen randomly and uniformly. Then $$p_0 < \left (1 -
\frac{s^2}{2^{k+2} k^3} \right)^q < e^{-s^2q/2^{k+2}k^3}$$

Since $$\frac{s^2q}{2^{k+2}k^3} = \Theta \left (\frac{n^2}{2^k k^3 q} 
\right) \rightarrow \infty \mbox{ as } k \rightarrow \infty ,$$ it follows that
$p_0$ approaches $0$ for large $k$. Thus the probability that some
arithmetic progression is monochromatic approaches $1$ as $k \rightarrow
\infty$. 
\end{proof}

\begin{theorem}
Let $N^{-}(k)= 2^{k/2} k^{1/2} g(k)$ where
$g(k) \rightarrow 0$ arbitrarily slowly as $k \rightarrow \infty$. Then,
the probability that a $2$-coloring of $\{1,2,\ldots,N^{-}(k)\}$ chosen
randomly and uniformly contains a monochromatic $k$-AP approaches $0$ as
$k \rightarrow \infty$. \\
\end{theorem}

\begin{proof}
Let $n=N^{-}(k)$, and let $E$ be the expected number
of monochromatic $k$-APs in a $2$-coloring of $\{1,2,\ldots,n\}$ chosen
randomly and uniformly. Note that there are $n^2(1+o(1))/(2k-2)$ $k$-APs
contained in $[1,n]$ and each of these is monochromatic with probability
$2^{1-k}$. By linearity of expectation, $E < k [g(k)]^2/(k-2)$. For $r \ge
0$, let $p_r$ be the probability that there are exactly $r$ monochromatic
$k$-APs in a random $2$-coloring. Then $E=p_1 + 2p_2 + 3p_3 + \ldots > 1 -
p_0$, so that $p_0 > (k-2-k[g(k)]^2)/(k-2)$. Thus, the probability that
some arithmetic progression is monochromatic approaches $0$ as $k
\rightarrow \infty$.  
\end{proof}

\section*{Acknowledgements}

The author thanks Tom Brown for bringing \cite{Ard05} to his attention,
and also for valuable comments on the first draft of this paper.

\end{document}